\documentclass{icmart}

\usepackage{amsrefs,bbm,verbatim}
\usepackage[T1]{fontenc}
\usepackage{pifont}
\usepackage{lipsum}
\usepackage{calligra}
\usepackage{chngcntr}
\counterwithout{footnote}{page}
\usepackage{pstricks, pst-node}
\psset{arrows=->, labelsep=3pt, mnode=circle}
\input pst-plot

\contact[asodin@princeton.edu]{Department of Mathematics, Princeton University, Fine Hall, Washington Road, Princeton, NJ 08544-1000, USA \\ \& \\ School of Mathematical Sciences, Tel Aviv University, Tel Aviv 6997801, Israel}

\numberwithin{equation}{section}
\numberwithin{figure}{section}

\newtheorem{theorem}{Theorem}[section]

\newtheorem{proposition}[theorem]{Proposition}

\newtheorem*{theorem*}{Theorem}

\theoremstyle{definition}

\newtheorem{remark}[theorem]{Remark}

\def\Wig{\sigma_{\mathrm{Wig}}}
\def\tils{\widetilde{s}}
\def\mlim{\mu_{\infty}}
\def\supp{\operatorname{supp}}
\def\Bern{\varepsilon}
\def\mes{\operatorname{mes}}
\def\sci{\text{\ding{33}}}
\def\rk{\operatorname{rank}}
\def\tr{\operatorname{tr}}

\def\Ai{\operatorname{Ai}}
\def\AD{\mathrm{A\!D}}
\def\xx{{\mathbf{x}}}
\def\XX{{\mathbf{X}}}
\def\xixi{{\boldsymbol{\xi}}}
\def\mm{{\mathbf{m}}}

\def\metr{{\mathbf{\rho}}}
\def\diag{\mathrm{diag}}
\def\offdiag{\mathrm{off\!\!-\!\!diag}}

\DeclareMathAlphabet{\mathcalligra}{T1}{calligra}{m}{n}

\title[Several applications of the moment method in random matrix theory]{Several applications of the moment method in random matrix theory}

\author[Sasha Sodin]{Sasha Sodin\footnote{Supported in part by NSF grant PHY-1305472}}

\begin{document}

\begin{abstract}
Several applications of the moment method in random
matrix theory, especially, to local eigenvalue statistics at the spectral edges,
are surveyed, with emphasis on a modification of the  method 
involving orthogonal polynomials.
\end{abstract}

\begin{classification}
60B20, 44A60
\end{classification}

\begin{keywords}
Moment method; random matrices; orthogonal polynomials
\end{keywords}

\maketitle

\bigskip
\hbox{{\vtop{\centering{\footnotesize\itshape \qquad\qquad To Yonatan Naim}}}}

\section{Introduction}

The goal of this article is to survey a few of the 
applications of the moment method (and its variants)
to the study of the spectral properties of random matrices, particularly, local eigenvalue statistics at
the spectral edges.

\vspace{1mm}\noindent
Section~\ref{s:moment} is a brief introduction  to the moment method, which we understand as the variety of ways to extract the properties of a measure $\mu$ from integrals of
the form
\begin{equation}\label{eq:moment} \int \xi^m \, d\mu(\xi)~. \end{equation}
Examples, selected from the narrow part of random matrix theory in which the author feels competent, are intended to illustrate two theses. First, the moment method can be  applied  beyond the  framework of weak convergence of a sequence of probability measures. Second, it is often convenient   to replace the monomials $\xi^m$ in (\ref{eq:moment}) with a better-conditioned sequence, such as the sequence of orthogonal polynomials with respect to a measure $\mlim$ which is an approximation to $\mu$. 

\vspace{1mm}\noindent
In Section~\ref{s:edges} we review some  applications   to the local eigenvalue statistics at the spectral edges,
starting from the work of Soshnikov \cite{Sosh}.
Tracy and Widom \cites{TW1,TW2} and Forrester \cite{For1} introduced the Airy point processes (see Section~\ref{s:edges.Wig})
and showed that they describe the limiting distribution  of the largest eigenvalues
for special families of large Hermitian random matrices with independent entries (the Gaussian invariant ensembles). 
Soshnikov \cite{Sosh} extended these results to Wigner matrices (Hermitian random matrices with independent entries
and no invariance assumptions).
In the terminology of Ibragimov and Linnik \cite{IL}*{Chapter VI}, the result of \cite{Sosh} 
is a limit theorem of collective character; it is one of the instances of the ubiquity (universality) of the Airy point processes
within and outside random matrix theory (as surveyed, for example, by Johansson \cite{Johdet}, Tracy and Widom \cite{TW3}, and 
Borodin and Gorin \cite{BorGor1}). 

In Section~\ref{s:edges.wigproc}, we consider  Wigner processes, a class of matrix-valued random processes. Informally, a random matrix $H(\xx)$ is attached to every point $\xx$ of an underlying space $\XX$. The statistical properties of the eigenvalues of every $H(\xx)$ are described by the theory of Wigner
matrices; the joint distribution of the eigenvalues of a tuple $(H(\xx_r))_{r=1}^k$ leads to limiting objects
which depend on the geometry of $\XX$ (which arises from the correlations of the matrix elements of $H$) in a non-trivial way. 
The moment method allows to derive limit theorems of collective character (such as Theorem~\ref{thm:wigproc}) pertaining to the
spectral edges of $H(\xx)$ (the result of \cite{Sosh}
corresponds to a singleton, $\# \XX = 1$).

In Section~\ref{s:edges.band}, we turn to the spectral 
edges of random band matrices. A
random band matrix is (\ref{eq:band}) a random
$N \times N$ Hermitian matrix with non-zero entries
in a band of width $W$ about the main diagonal. When
$W$ is small, a band matrix inherits the structure
of the integer lattice $\mathbb{Z}$; when $W$ is large,
it is similar to a Wigner matrix. The threshold at which the local eigenvalue statistics in the bulk of the spectrum exhibit a crossover is
described by precise conjectures (see Fyodorov and Mirlin \cites{FM,FM2}, Spencer \cites{Sp-banded,Sp-SUSY}). 
The moment method allows to prove the counterpart of these conjectures for the spectral edges (the result of \cite{Sosh} corresponds to the special case $W = N$).

\vspace{1mm}\noindent
The content of Section~\ref{s:moment} is mostly  known. The modified moment method of Section~\ref{s:modif} is a version of self-energy renormalisation in perturbation theory (see Spencer \cite{Sp-banded}), related to the arguments of Bai and Yin \cite{BY}. Orthogonal polynomials were explicitly used in this context in the work of Li and Sol\'e \cite{LiSole}, and further in \cite{meop} (where more references may be found). Some  observations are incorporated from \cite{FS1}. The content of Section~\ref{s:edges.wigproc}  is an extension of \cite{mecorner}, whereas  Section~\ref{s:edges.band} is based on \cite{meband1}. The proofs of the results stated in both of these sections build on the combinatorial arguments of \cite{FS1}.

\section{Preliminaries and generalities}\label{s:moment} The moment method is the collection of techniques inferring 
the properties of a measure $\mu$ on the $k$-dimensional space $\mathbb{R}^k$ from the moments
\begin{equation}
s(m_1, \cdots,m_k; \mu) = \int_{\mathbb{R}^k} 
\xi_1^{m_1} \cdots \xi_k^{m_k} \, d\mu(\xi)
\quad (m_1, \cdots, m_k =0,1,2,\cdots)~.
\end{equation}
Introduced by Chebyshev as a means to establish  Gaussian approximation for 
the distribution of a sum of independent random
variables, the moment method achieved its first major success with the
proof, given by Markov \cite{Markov}, of Lyapunov's Central Limit Theorem and its extension to 
sums of weakly dependent random variables. Some of the more recent applications are surveyed by
Diaconis \cite{Dia}.

\subsection{Convergence of probability measures}\label{s:convprob} In the traditional
setting of the moment method, one 
considers a sequence of probability measures $(\mu_N)_{N \geq 1}$ on $\mathbb{R}^k$. Suppose
that the limit
\begin{equation}\label{eq:momentconv}
s(m_1, \cdots, m_k) = \lim_{N \to \infty} s(m_1, \cdots, m_k; \mu_N)
\end{equation}
exists for every $m_1, \cdots, m_k \geq 0$. Then
the sequence $(\mu_N)_{N \geq 1}$ is tight, i.e.\ 
precompact in weak topology (defined by bounded continuous functions), and every one of its limit
points $\mu$ satisfies
\begin{equation}\label{eq:moments}  s(m_1, \cdots, m_k; \mu) = s(m_1, \cdots, m_k)
\quad (m_1, \cdots,m_k \geq 0)~. 
\end{equation}
If, for example,
\begin{equation}\label{eq:expm}
s(m_1, \cdots, m_k) \leq \prod_{r=1}^k (C m_r)^{m_r} \quad (m_1, \cdots, m_k \geq 0)~,
\end{equation}
the moment problem (\ref{eq:moments}) is determinate, i.e.\ there is a unique measure $\mlim$
on $\mathbb{R}^k$ satisfying (\ref{eq:moments}); in
this case the convergence of moments (\ref{eq:momentconv}) implies that
$\mu_N \to \mlim$ in weak topology; cf.\ Feller 
\cite{Fel}*{\S VIII.6}.

Hardy's sufficient condition (\ref{eq:expm}) may be somewhat
relaxed; we refer to the addenda to the second 
chapter of the book \cite{Akh} of Akhiezer for 
various sufficient criteria for determinacy in the case
$k=1$, and to the survey of Berg \cite{Berg} for
some extensions to $k \geq 1$.

\subsubsection{Random measures}\label{sub:randmeas}

Suppose $(\mu_N)_{N \geq 1}$ is a sequence of random measures on $\mathbb{R}^\ell$ (i.e.\ random variables
taking values in the space of Borel probability measures). We denote $\xixi = (\xi_1,  \cdots, \xi_\ell)$, $\mm = (m_1, \cdots, m_\ell)$,
and $\xixi^\mm = \xi_1^{m_1} \cdots \xi_\ell^{m_\ell}$; thus
 $\xixi_r^{\mm_r} =  \xi_{r,1}^{m_{r,1}} \cdots \xi_{r,\ell}^{m_{r,\ell}}$. If 
\begin{equation}\label{eq:randmconv}\mathbb{E} \int_{\mathbb{R}^{\ell k}} \prod_{r=1}^k \left[ d\mu_N(\xixi_r) \xixi_r^{\mm_r}  \right] \to
\mathbb{E} \int_{\mathbb{R}^{\ell k }} \prod_{r=1}^k  \left[ d\mlim(\xixi_r)  \xixi_r^{\mm_r}  \right]  \quad (k \geq 1)
\end{equation}
for some random measure $\mlim$ on $\mathbb{R}^\ell$, and the moment problem
for every moment measure $\mathbb{E} \mlim^{\otimes k}$ is determinate, then
\begin{equation}\label{eq:convmm}
\mathbb{E} \mu_N^{\otimes k} \overset{\text{weak}}{\longrightarrow}
\mathbb{E} \mlim^{\otimes k} \quad (k \geq 1)~.
\end{equation}
If, for every Borel set $K \in \mathcal{B}(\mathbb{R}^\ell)$, the
moment problem for the distribution of $\mlim(K)$
is determinate, then (\ref{eq:convmm}) implies
that $\mu_N \to \mlim$ (weakly in distribution).
Finally, if  $\mlim$ is deterministic (i.e.\ its distribution is supported on one deterministic measure), it is sufficient to verify
(\ref{eq:convmm}) for $k=1,2$. 
See further Zessin \cite{Zes}.

\subsection{Example: Wigner's law}\label{sub:Wigner}
The application (going back to Chebyshev) of the moment method to sums of independent random variables is based on the identity
\begin{equation}\label{eq:idclt} \left[ \sum_{j=1}^N X_j\right]^m = 
\sum_{m_1 + \cdots + m_N = m} \frac{m!}{m_1! \cdots m_N!} \prod_{j=1}^N  X_j^{m_j}
\end{equation}
expressing powers of a sum of numbers as a sum over partitions.
 Similarly,
the application (going back to Wigner) of the moment method to random
matrix theory is based on the relation
\[\begin{split} \tr H^m
&=\sum_{p = (u_0, u_1, \cdots, u_{m-1}, u_m = u_0)} \prod_{j=0}^{m-1} H(u_j, u_{j+1}) \\
&= \sum_p \prod_{1 \leq u \leq v \leq N} 
H(u, v)^{\# \{ (u_j, u_{j+1}) = (u, v)\}}
\overline{H(u, v)}^{\# \{ (u_j, u_{j+1}) = (v, u)\}}
\end{split}\]
expressing traces of powers of an  Hermitian  matrix $H = (H(u,v))_{u,v=1}^N$ 
as a sum over paths.

\vspace{1mm}\noindent
Let $\big(G_N = (V_N, E_N)\big)_{N \geq 1}$ be a sequence of graphs, 
so that $G_N$ is $\kappa_N$-regular (meaning that every
vertex is adjacent to exactly $\kappa_N$ edges),
and  the connectivity $\kappa_N$ tends
to infinity:
\begin{equation}\label{eq:largekappa'} \lim_{N \to \infty} \#V_N = \infty, \quad  \lim_{N \to \infty} \kappa_N = \infty~.
\end{equation}
For every $N$, consider a random matrix $H^{(N)} = (H(u,v))_{u, v \in V_N}$ with rows
and columns indexed by the elements of $V_N$,  so that
$(H(u, v))_{u,v \in V_N}$  are independent up to the constraint
$H(v, u) = \overline{H(u, v)}$;
the diagonal entries $\{H(u, u)\}_{u}$ are sampled from a distribution
$\mathcal{L}_\text{diag}$ on $\mathbb{R}$ satisfying
\begin{equation}\label{eq:diag} \mathbb{E} H(u, u) = 0~, \quad \mathbb{E} H(u, u)^2 < \infty~;
\end{equation}
 the off-diagonal entries $\{H(u, v)\}_{(u,v) \in E_N}$
are sampled from a distribution $\mathcal{L}_\text{off-diag}$ on $\mathbb{C}$ satisfying
\begin{equation}\label{eq:offdiag} \mathbb{E} H(u, v) = 0~, \quad \mathbb{E} |H(u, v)|^2 = 1~;
\end{equation} 
and all the other entries $H(u, v)$ are set to zero.

Let $\xi_1^{(N)} \geq \xi_2^{(N)} \geq \cdots \xi_{\# V_N}^{(N)}$
be the eigenvalues of $H^{(N)}$,
and let
\[ \mu_N = \frac{1}{\# V_N} \sum_{j=1}^{\# V_N}  \delta\left(\xi - \frac{\xi_j^{(N)}}{2\sqrt{\kappa_N-1}}\right)~. \]
(The scaling is natural since, for instance, the $\ell_2$ norm of every column of the $N \times N$ matrix is of order $\sqrt{\kappa_N}$.)

\begin{theorem*}[Wigner's law]  
In the setting of this paragraph (i.e.\ assuming (\ref{eq:largekappa'}), (\ref{eq:diag}), (\ref{eq:offdiag})), the sequence of random measures $(\mu_N)_N$
converges (weakly, in distribution) as $N \to \infty$ to
the (deterministic) semicircle measure $\Wig$ with density
\begin{equation}\label{eq:Wig}\frac{d\Wig}{d\xi} = \frac{2}{\pi} \sqrt{(1-\xi^2)_+}~. \end{equation}
\end{theorem*}

Wigner considered \cites{Wig1,Wig2,Wig3} the case when $G_N$ is the complete
graph on $N$ vertices (Wigner matrices), and the entries satisfy some
additional assumptions, the important of them being that all moments are  finite. Wigner's argument is based on the moment method.

Bogachev, Molchanov, and Pastur \cite{BMP} observed (in the context of random band matrices) that a similar argument can be applied as long as (\ref{eq:largekappa'}) is
satisfied. The first argument for Wigner matrices
without additional restrictions on the distribution
of the entries was given by Pastur \cite{P1}, using the
Stieltjes transform method introduced by Marchenko
and Pastur \cites{MP,MP2} (see Pastur \cite{P2} and the book of Pastur and Shcherbina \cite{PSbook} for  some of the further applications of the method).  Khorunzhiy,
Molchanov, and Pastur \cite{KMP} applied the Stieltjes transform method
to prove Wigner's law for random band matrices; their argument is applicable in the
 setting described here.

Let us outline a proof of Wigner's law in the form stated above,
following  \cite{BMP}
(and incorporating Markov's truncation argument \cite{Markov}). We refer for details to the book of Anderson, Guionnet, and Zeitouni \cite{AGZ}*{Chapter 2.1},
where similar arguments are also applied to  questions such as the Central Limit Theorem for 
linear statistics $\phi(\xi_1^{(N)}) + \cdots + \phi(\xi_{\# V_N}^{(N)})$.

\begin{proof}[Proof of Wigner's law]
Due to (\ref{eq:largekappa'}), (\ref{eq:diag}) and (\ref{eq:offdiag}) one can find
a sequence $\delta_N \to +0$ so that 
\[ \mathbb{E} |H(u, v)|^2 \mathbbm{1}_{|H(u,v)|\geq \delta_N \sqrt{\kappa_N}} \leq \delta_N~.\]
Consider the  matrix $H_\sci^{(N)}$ with truncated matrix elements
\[ H_\sci(u, v) = 
\begin{cases}
H(u, v)~, & |H(u, v)| \leq \sqrt{\kappa_N} \\
0~, &|H(u, v)| > \sqrt{\kappa_N}
\end{cases}~.
 \]
Then
\[ \mathbb{P} \left\{ H_\sci(u, v) \neq H(u, v) \right\} \leq \delta_N \kappa_N^{-1}~, \]
whence, bounding rank by the number of non-zero matrix elements and applying the Chebyshev inequality,
\[ \mathbb{P} \left\{ \rk (H_\sci^{(N)} - H^{(N)}) \geq \delta_N^{1/2} \# V_N \right\}
\leq \delta_N^{1/2}~.\]
For any  $\xi \in \mathbb{R}$, the interlacing property of rank-one perturbation yields
\begin{equation}\label{eq:rk1}
\big| \mu_N(-\infty,\xi] - \mu_{\sci,N}(-\infty,\xi]\big| \leq \#V_N^{-1} \rk(H_\sci^{(N)} - H^{(N)})~,
\end{equation}
therefore it is sufficient to establish the
result for $H_\sci^{(N)}$ in place of $H^{(N)}$. For large $N$, the elements of $H_\sci^{(N)}$ enjoy the following estimates:
\begin{align}
\label{eq:trest1}\left|E H_\sci(u, v)\right| &\leq \delta_N \kappa_N^{-\frac12}~; \\
\label{eq:trest2}\left|\mathbb{E} |H_\sci(u, v)|^2 - 1 \right| &\leq \delta_N \quad ((u,v) \in E_N)~; 
	\quad \mathbb{E} |H_\sci(u, u)|^2 \leq \mathrm{const}~;\\
\label{eq:trest3}\mathbb{E} \left| H_\sci(u, v) \right|^k &\leq 2 \delta_N \kappa_N^{\frac{k-2}{2}} \quad
(k \geq 3)~.
\end{align}

Next, consider the expansion
\begin{equation}\label{eq:trm}
\begin{split}
s(m_1, \cdots, m_k; \mathbb{E} \mu_{\sci,N}^{\otimes k})
&= \mathbb{E} \prod_{r=1}^k \int \xi^{m_r} d\mu_{\sci,N}(\xi) \\
&= \mathbb{E} \frac{1}{(\# V_N)^k} \prod_{r=1}^k
\tr \left( \frac{H_\sci^{(N)}}{2 \sqrt{\kappa_N - 1}} \right)^{m_r} \\
&= \sum \frac{1}{(\# V_N)^k} \mathbb{E} \prod_{r=1}^k \prod_{j=0}^{m_r - 1} 
   \frac{H_\sci(u_{r,j}, u_{r ,j+1})}{2\sqrt{\kappa_N - 1}}~,
\end{split}
 \end{equation}
where the sum is over $k$-tuples of closed paths 
\[\begin{split}\begin{cases}
 &u_{1,0}, u_{1, 1}, \cdots, u_{1, m_1-1}, u_{1, m_1} \\
&u_{2,0}, u_{2, 1}, \cdots, u_{2, m_2-1}, u_{2, m_2} \\
&\cdots\\
&u_{k,0}, u_{k, 1}, \cdots, u_{k, m_k-1}, u_{k,m_k} \\
  \end{cases} \quad
 \left[  u_{1, m_1} = u_{1,0}, \cdots, u_{k, m_k} = u_{k,0} \right] \end{split}\]
in the augmented (multi-)graph $G_N^+ = (V_N, E_N^+)$, 
$E_N^+ = E_N \bigcup \left\{ (u, u) \, \mid \, u \in V_N \right\}$.
Two such $k$-tuples are called isomorphic if one is obtained 
from one another by a permutation of the vertices $V_N$. For example, the pair $(1 \, 3 \, 1~,\, 2 \, 1 \, 4 \, 2)$ is isomorphic to $(7 \, 4 \, 7~,\,1 \, 7 \, 2 \, 1)$.

According to (\ref{eq:trest1}), (\ref{eq:trest2}) and (\ref{eq:trest3}),
the contribution of an isomorphism class consisting of $k$-tuples spanning 
a graph $\mathfrak{g}$ with $\mathfrak{v}$ vertices and $\mathfrak{e}$ edges, 
of which $\mathfrak{e}_2$ are traversed exactly twice, is bounded by $\kappa_N^{\mathfrak{v} - k} (\delta_N)^{\mathfrak{e} - \mathfrak{e}_2} 
(\kappa_N/\mathrm{const})^{-\mathfrak{e}}$.

The graph $\mathfrak{g}$ has at most $k$ connected components, whence
\begin{equation}\label{eq:Euler} \mathfrak{v} -\mathfrak{e} \leq k~,
\end{equation}
with equality for graphs which are vertex-disjoint unions of $k$ trees. 
For fixed $m_1, \cdots, m_k$, the number of isomorphism classes remains bounded as $N \to \infty$, 
therefore the limit
of (\ref{eq:trm}) is given by the contribution of vertex-disjoint $k$-tuples of paths 
corresponding to graphs with
\begin{equation}\label{eq:leading}
 \mathfrak{v} - \mathfrak{e} = k~, \quad \mathfrak{e}_2 = \mathfrak{e}~. 
\end{equation}
 
\begin{figure}[ht]
\vspace{2.2cm}
\setlength{\unitlength}{1cm}
\begin{pspicture}(-1,0)
\psset{arrows=-}
\rput(1, .9){\small 1}
\psline(1.2,1)(2, 1)
\rput(1.9, 1.2){\small 2}
\psline(2,1)(2.5,2)
\psline(2.5,2)(2.6,2)
\rput(2.5,2.2){\small 3}
\psline(2.6, 2)(2.1,1)
\psline(2.1, 1)(4,1)
\rput(3.9,1.2){\small 5}
\psline(4,1)(4.5, 2)
\psline(4.5,2)(4.6,2)
\rput(4.5,2.2){\small 6}
\psline(4.6,2)(4.1,1)
\psline(4.1, 1)(4.6, 0)
\psline(4.6,0)(4.5,0)
\rput(4.7, .2){\small 7}
\psline(4.5,0) (4,.9)
\psline(4,.9)(3,.9)
\psline(3,.9)(3.5, -.1)
\rput(3.6,.1){\small 8}
\psline(3.5, -.1)(4, -.1)
\psline(4,-.1)(4,-.2)
\rput(3.9,.1){\small 9}
\psline(4, -.2)(2.4, -.2)
\psline(2.4,-.2)(2.4,-.1)
\rput(2.4,.1){\small 10}
\psline(2.4, -.1)(3.4, -.1)
\psline(3.4, -.1)(2.9,.9)
\rput(2.9, 1.2){\small 4}
\psline(2.9,.9)(1.2,.9) %
\psline(6.2,1)(8, 1)
\rput(6,.9){\small 1}
\rput(6.7,1.2){\small 2}
\rput(7.4,1.2){\small 3}
\rput(8.3,.9){\small 4}
\psarcn(8.8,1){.8}{180}{195}
\rput(8.5,1.4){\small 5}
\rput(9.1,1.4){\small 6}
\rput(9.3,.9){\small 7}
\rput(8.75,.5){\small 8}
\psline(8,.8)(8.1,.92)
\psarcn(8.8,1){.7}{185}{195}
\psline(8,.92)(8.1,.8)
\psline(8,.9)(6.2,.9) %
\end{pspicture}
\caption{The tree-like path $1\, 2 \, 3 \, 2 \, 4 \, 5 \, 6 \, 5 \, 7 \, 5 \,
4 \, 8 \, 9 \, 8 \, 10 \, 8 \, 4 \, 2 \, 1$ with $\mathfrak{v}=10$ and
$\mathfrak{e}=\mathfrak{e}_2=9$ (left) and the non-backtracking path
$1 \, 2 \, 3 \, 4 \, 5 \, 6 \, 7 \, 8 \, 4 \, 5 \, 6 \, 7 \, 8 \, 4 \, 3 \, 2 \, 1$
with $\mathfrak{v} = 8$ and $\mathfrak{e}=\mathfrak{e}_2=8$ (right). Among the two, only the first one contributes to the semi-circle limit.}\label{fig:forest}
\vspace{2mm}
\end{figure}
Every path in such a $k$-tuple is tree-like (see Figure~\ref{fig:forest}, left); each isomorphism class contributes
$2^{-\sum_p m_p}$ (due to (\ref{eq:trest2})), and the number of classes is given by a product
of Catalan numbers:
\[ \prod_{r=1}^k \begin{cases}
 \frac{2}{m_r + 2} \binom{m_r}{m_r/2}~, &\text{$m_r$ is even}\\
0~, &\text{$m_p$ is odd}                 
 \end{cases} = \prod_{r=1}^k \big[ 2^{m_r} s(m_r; \Wig)\big]~.\]
Thus 
\begin{equation}\label{eq:wigconv}
\lim_{N \to \infty} s(m_1, \cdots, m_k; \mathbb{E} \mu_{\sci, N}^{\otimes k}) =
s(m_1, \cdots, m_k; \Wig^{\otimes k})~.
\end{equation}
Applying the relation (\ref{eq:wigconv}) with $k=1,2$, we conclude (cf.\ Section~\ref{sub:randmeas}) that $\mu_{\sci,N}$ converge to $\Wig$
weakly in distribution, and thus (by (\ref{eq:rk1})) so do $\mu_{N}$.
\end{proof}

\subsection{Some quantitative aspects}
Whenever the moment convergence
(\ref{eq:momentconv}) is a consequence of the stronger property
\begin{equation}\label{eq:momentconv+}
s(m_1, \cdots, m_k; \mu_N) =s(m_1, \cdots, m_k; \mlim)
\quad (N \geq N_0(m_1, \cdots, m_k))~,
\end{equation}
the arguments quoted in Section~\ref{s:convprob} can be recast in quantitative form. 
This is illustrated by the following inequality due to  Sonin \cite{Son}. Let $\gamma$ be the Gaussian
measure,
\[ \frac{d\gamma}{d\xi} = \frac{1}{\sqrt{2\pi}} \exp(-\xi^2/2) \quad (\xi \in \mathbb{R})~,\]
and assume that
\begin{equation}\label{eq:matchgauss}
 s(m; \mu_N) = s(m; \gamma)
\left[ = \begin{cases} 0, &\text{$m$ is odd} \\
\frac{m!}{(\frac{m}{2})!2^{\frac{m}{2}}}, &\text{$m$ is even}
\end{cases}\right]
\quad (N \geq N_0(m))~.
\end{equation}
Then 
\begin{equation}\label{eq:sonin} \sup_{\xi \in \mathbb{R}} \left| \mu_N(\xi) -
\gamma(\xi) \right| \leq \sqrt{\frac{\pi}{m-1}}
\quad (N \geq \max_{m' \leq m} N_0(m'))~.\end{equation}
Measures $\mu_N$ of random matrix origin for which (\ref{eq:matchgauss}) holds
may be found in the survey of Diaconis \cite{Dia}. Inequalities of the form
(\ref{eq:sonin}) may be also derived for other measures $\mlim$ (see Akhiezer \cite{Akh}*{Section~II.5.4}
for the general framework of  Chebyshev--Markov--Stieltjes inequalities,
and Krawtchouk \cite{Kraw} for additional examples).

Similar inequalities can be derived for $k > 1$. On
the other hand, already in the setting of the Central
Limit Theorem for sums of independent random variables, (\ref{eq:matchgauss}) is not valid (unless the
addends are Gaussian themselves); the correct relation
$ s(m; \mu_N) \approx s(m; \gamma)$, 
even with the optimal dependence of the error term on $m$ and $N$, yields a 
poor bound on the rate of convergence of $\mu_N$ to $\gamma$ (the sharp Berry --Esseen 
bound, see Feller \cite{Fel}*{\S XVI.5}, was  proved by the Fourier-analytic approach). The reason is that monomials
form an ill-conditioned basis; see Gautschi \cite{Gau} for a discussion of computational aspects (and of remedies similar to the one discussed in the
next section).

\subsection{A modification of the moment method}\label{s:modif}
The following modification makes the moment method
better conditioned. Let $(\mu_N)_{N \geq 1}$ be a sequence of probability measures on $\mathbb{R}$, and suppose $\mlim$ is a
candidate for the weak limit of the sequence  $(\mu_N)_{N \geq 1}$. Let  $P_n(\xi)$ ($n=0,1,2,\cdots$)
be the orthogonal polynomials with respect to 
$\mlim$: 
\[ \deg P_n = n~, \quad \int P_n(\xi) P_{n'}(\xi) d\mu_\infty(\xi) = \delta_{nn'}~.\] 
Also set
\begin{equation}\label{eq:tils} \tils(n; \mu; \mlim) = \int_{-\infty}^\infty 
P_n(\xi) d\mu(\xi)~.\end{equation}
Then the convergence of moments \begin{equation}\label{eq:momentconv'}
\lim_{N \to \infty} s(m; \mu_N)
= s(m; \mlim) \quad (m \geq 0)
\end{equation}
is equivalent to 
\begin{equation}\label{eq:pconv}
\lim_{N \to \infty} \tils(n; \mu_N; \mlim)
= \delta_{n0} \quad (n \geq 0)~.
\end{equation}
Thus (\ref{eq:pconv}) implies that $\mu_N \to \mlim$,
provided that the moment problem for $\mlim$ is determinate.

While the modification of the moment
method advertised here seems to have no general counterpart in dimension $k>1$,
in  the special case when
$\mlim$ is the \mbox{$k$-th}  power of a one-dimensional measure with orthogonal polynomials $P_n$
we define:
\[ \tils(n_1, \cdots, n_k; \mu; \mlim)
= \int_{\mathbb{R}^k} \prod_{r=1}^k P_{n_r}(\xi_r) \,
d\mu(\xi)~. \]

\subsubsection{A random matrix example}\label{subsub:mod}

If  $X_1, \cdots, X_n$ are independent random variables with zero mean, unit variance, and finite moments, one may give a combinatorial interpretation to
\begin{equation}\label{eq:idclt1} \mathbb{E} \frac{1}{\sqrt{n!}} 
\operatorname{He}_n \left[ \frac{X_1 + \cdots + X_N}{\sqrt{N}}\right]~,
\end{equation}
where 
\[ \operatorname{He}_n(\xi) = (-1)^n e^{\xi^2/2} 
\frac{d^n}{d\xi^n} e^{-\xi^2/2} \]
are the Hermite polynomials; the three-term recurrent relation
\[ \operatorname{He}_{n+1}(\xi) = \xi \operatorname{He}_n (\xi) - n \operatorname{He}_{n-1}(\xi)\]
eliminates  the asymptotically leading terms of the moments 
(\ref{eq:idclt}) of $X_1 + \cdots + X_N$. Here we focus on a different example,
pertaining to random matrices of the form considered in Section~\ref{sub:Wigner}.

Denote
\[ P_n^{(\kappa)}(\xi) = U_n(\xi) - \frac{1}{\kappa - 1} U_{n-2}(\xi)~, \]
where 
\[ U_n(\cos \theta) = \frac{\sin ((n+1)\theta)}{\sin \theta} \]
are the Chebyshev polynomials of the second kind (orthogonal with respect to $\Wig$), and  $U_{-1} \equiv U_{-2} \equiv 0$. 
Let $G = (V, E)$ be a regular graph of connectivity $\kappa$, and let $H$ be an $\#V \times \#V$ Hermitian matrix, 
such that
\begin{equation}\label{eq:unimod}
|H(u, v)| = \mathbbm{1}_{(u,v) \in E}~, \quad (u, v \in V)~.
\end{equation}
The three-term recurrent relation 
\[P_{n+1}^{(\kappa)} (\xi)= 2\xi P_n^{(\kappa)}(\xi) 
- (1 + (\kappa-1)^{-1}\mathbbm{1}_{n=1}) P_{n-1}^{(\kappa)}(\xi)\]
for $P_n^{(\kappa)}$ leads to
\begin{proposition}[cf.\ \cite{meop}*{Lemma~2.7}, \cite{FS1}*{Claim~II.1.2}]\label{prop:cl}
For any Hermitian matrix $H$ satisfying  (\ref{eq:unimod}), 
\begin{equation}\label{eq:main}
P_n^{(\kappa)}\left[\frac{H}{2\sqrt{\kappa-1}}\right](u, v)
= \sum \prod_{j=1}^n \frac{H(u_j, u_{j+1})}{\sqrt{\kappa-1}}~,
\end{equation}
where the sum is over paths $u_0, u_1,\cdots, u_{n-1}, u_n$ in $G$ from $u_0 = u$ to $u_n = v$
which satisfy the non-backtracking condition $u_j \neq u_{j+2}$ ($0 \leq j \leq n-2$).
\end{proposition}
Consider a sequence of random matrices $H^{(N)}$ associated to a sequence of graphs $G_N$ with $\kappa_N \to \infty$ as in Section~\ref{sub:Wigner}; let us assume that the
entries of $H$ satisfy the unimodality assumptions
(\ref{eq:unimod}). A non-backtracking path can not be tree-like (see Figure~\ref{fig:forest}), therefore the modified
moments tend to zero; this provides an alternative proof to
Wigner's law in the form of  Section~\ref{sub:Wigner} under the additional assumptions
(\ref{eq:unimod}).

\vspace{1mm}\noindent
The generalisation of Propostion~\ref{prop:cl} to
matrices which do not satisfy (\ref{eq:unimod}) is somewhat technical, and we do not present it here. In the context
of Wigner (and sample covariance) matrices, it
is described in \cite{FS1}*{Part~III}; for the (more involved) case of band matrices we refer to the 
work of Erd\H{o}s and Knowles \cite{EK2}.

\subsubsection{Advantages of modified moments}\label{s:adv}

Although the convergence of modified moments (\ref{eq:pconv}) is equivalent to
the convergence of moments (\ref{eq:momentconv'}), quantitative forms of the
former yield better estimates on the rate of convergence $\mu_N \to \mlim$.
As an illustration, we recall 
a variant of the Erd\H{o}s--Tur\'an inequality 
\cite{ET1} proved in \cite{FS2}.
Consider again the semi-circle measure $\Wig$ with
density (\ref{eq:Wig}).
\begin{proposition}[{\cite{FS2}*{Proposition~5}}]\label{p:et} Let $\mu$ be a probability
measure on $\mathbb{R}$. Then, for any $\xi \in \mathbb{R}$ and any $n_0 \geq 1$, 
\[ \left| \mu(-\infty, \xi] - \Wig(-\infty,\xi] \right| \leq C \left\{ \frac{\rho(\xi; n_0)}{n_0}
+ \sqrt{\rho(\xi; n_0)} \sum_{n=1}^{n_0} \frac{|\tils(n; \mu; \Wig)|}{n}\right\}~,\]
where $C>0$ is a numerical constant, and 
$\rho(\xi; n_0) = \max(1-|\xi|, n_0^{-2})$.
\end{proposition}
The original Erd\H{o}s--Tur\'an inequality provides
a bound of similar structure for the measure with
density 
${d\mlim}/{d\xi} = \pi^{-1}\left((1-\xi^2)_+\right)^{-1/2}$
(in this case, $\rho(\xi; n_0)$ should be replaced with $1$.)\footnote{A similar inequality for the Gaussian measure, combined with a careful estimate
of the modified moments (\ref{eq:idclt1}), could perhaps yield 
a proof of the Berry--Esseen theorem along the
lines suggested by Chebyshev.}

\subsection{Convergence of rescaled probability measures}\label{s:resc}

The rescaling $R_{\eta}^{\xi_0} [\mu]$ of a measure $\mu$ on $\mathbb{R}^k$ about $\xi_0 \in \mathbb{R}^k$
by $\eta > 0$ is defined by
\begin{equation}\label{eq:rescop} R_{\eta}^{\xi_0} [\mu](K) = \mu(\eta(K - \xi_0)) \quad \big(K \in \mathcal{B}(\mathbb{R}^k)\big)~. 
\end{equation}
In a class of questions  outside the narrow framework of Section~\ref{s:convprob}, one 
is interested in vague limits (weak limits with respect
to the topology defined by compactly supported continuous functions) of 
\begin{equation}\label{eq:resc} \Big(\epsilon_N^{-1} R_{\eta_N}^{\xi_0}[\mu_N]\Big)_{N \geq 1}~, \end{equation}
where $(\mu_N)_{N \geq 1}$  is a sequence of probability measures on $\mathbb{R}^k$, $\xi_0 \in \mathbb{R}^k$, and  two sequences 
$\epsilon_N, \eta_N \to +0$ determine the scaling of $\mu_N$ on the  value ($\updownarrow$) and 
variable ($\leftrightarrow$)
axes, respectively.

\subsubsection{Edges (corners) of the support}\label{subsub:corn} Moments allow to study the rescaling of $\mu_N$ about a point $\xi_0$ which is close to the corners of the cube supporting $\mu_N$.  Variants of this observation were
used, for example, by Sinai and Soshnikov \cites{SinSosh1,SinSosh2}. 

Assume that we are given a sequence $(\mu_N)_{N \geq 1}$ 
of probability measures on $\mathbb{R}^k$,  two sequences $\epsilon_N, \eta_N \to + 0$ which determine the scaling (\ref{eq:resc}),
and $2^k$ continuous functions  
$\phi_\Bern: (\alpha_0,\infty)^k \to \mathbb{R}_+$  ($\Bern \in \{-1,1\}^k$) which will describe the
limiting Laplace transform at the $2^k$ corners of the cube.

\begin{proposition}\label{p:lapl}
Suppose
\[ \epsilon_N^{-k} s(m_{1,N}, \cdots, m_{k,N}; \mu_N) - \sum_{\Bern \in \{-1, 1\}^k} \prod_{r=1}^k \Bern_r^{m_r} \, \phi_\Bern(\alpha_1, \cdots,
\alpha_k) \longrightarrow 0 
\quad (N \to \infty) \]
for any sequence $(m_{1, N},\cdots,m_{k,N})_{N \geq 1}$ for which 
\[ \lim_{N \to \infty}\eta_N m_{r,N} =  \alpha_r > \alpha_0 \quad(1 \leq r \leq k).\]
Then, for any $\Bern \in \{-1,1\}^k$, the sequence
$(\epsilon_N^{-1}R_{\eta_N}^\Bern[\mu_N])_{N \geq 1}$
converges vaguely to a measure $\nu^\Bern$
which is uniquely determined by the equations
\[ \int \exp(\alpha_1 \lambda_1+ \cdots + \alpha_k \lambda_k) d\nu^\Bern(\Bern_1 \lambda_1, \cdots,
\Bern_k \lambda_k) = \phi_\Bern(\alpha_1, \cdots, \alpha_k)
\quad (\alpha \in (\alpha_0, \infty)^k)~.\]
\end{proposition}
\begin{remark}\label{rem:top} Convergence actually holds in the stronger topology defined by continuous functions
supported (for some $R > 0$) in
\[ \prod_{r=1}^k 
\begin{cases} 
(-R, \infty)~, & \Bern_r = 1 \\
(-\infty, R)~, & \Bern_r = -1
\end{cases}~. \]
\end{remark}

The counterparts of Proposition~\ref{p:lapl} for modified moments depend
on the structure of the limiting measure $\mlim$. For the case $\mlim = \Wig^{\otimes k}$ such a statement was proved
in \cite{meband1}*{Section~6}. It is somewhat technical, and we do not reproduce it here; instead of the Laplace transform, the
limiting measures $\nu^\Bern$ are characterised
in terms of the transform
\begin{equation}\label{eq:transf} \int_{\mathbb{R}^k} \prod_{r=1}^k 
\frac{\sin \alpha_r \sqrt{-\lambda_r}}{ \sqrt{-\lambda_r}} \, d\nu^\Bern(\Bern_1 \lambda_1,
\cdots, \Bern_k \lambda_k)
\end{equation}
(which becomes convergent after a certain regularisation).  The system of functions
$\lambda \mapsto \frac{\sin \alpha \sqrt{-\lambda}}{\sqrt{-\lambda}}$ 
forms a continuous analogue of orthogonal polynomials (as introduced by Krein,
see Denisov \cite{Den}) with respect to the measure $\frac{2\sqrt{2}}{\pi} \sqrt{-\lambda_-}$ (obtained by
rescaling $\Wig$ about $\xi_0=1$). 

Uniqueness theorems for the
transform (\ref{eq:transf}) were proved (in dimension $k=1$) in the
1950-s by Levitan \cite{Lev}, Levitan and Meiman \cite{LevMeim}, and Vul \cite{Vul} (listed in order of increasing
generality);
the argument in \cite{meband1} builds on \cite{Lev}.

One advantage of the approach based on modified moments is that, for a measure supported on several intervals, it allows to consider the rescaling about  edges (corners) which are not maximally distant from the origin, and even internal edges. In the context of random matrices, this was exploited in \cite{FS1}.

\subsubsection{Interior points of the support}\label{s:resc.int}

If $\xi_0$ is  an interior point of the support of $\mlim$, it seems impossible to extract any information regarding the measures $\epsilon_N^{-1} R_{\eta_N}^{\xi_0}[\mu_N]$
from the asymptotics of the moments of $\mu_N$. The modified moments $\tils$
carry such information. For example, Proposition~\ref{p:et} shows that if one
can find a sequence $(n_0(N))_{N \geq 1}$ so that 
\begin{equation}\label{eq:etcond}
\lim_{N \to \infty} \epsilon_N n_0(N) = +\infty~,\quad
\lim_{N \to \infty} n_0(N) \sum_{n = 1}^{n_0(N)} \frac{|\tils(n; \mu_N; \Wig)|}{n} = 0~,
\end{equation}
then
\begin{equation}\label{eq:tomes}
 \epsilon_N^{-1} R_{\epsilon_N}^{\xi_0} [\mu_N] \overset{\text{vague}}{\underset{N \to \infty}{\longrightarrow}} 
\frac{1}{\pi} \sqrt{1-\xi_0^2} \,\, \mes \qquad (-1 < \xi_0 < 1)
\end{equation}
(where $\mes$ is the Lebesgue measure on the real line). 

Let us briefly comment on the shorter scales $\epsilon_N$, for which (\ref{eq:etcond})
fails. The challenge is to give meaning to
the expansion
\begin{equation}\label{eq:L2exp}
\mu_N[\xi', \xi''] \sim \sum_{n \geq 0} \tils(n; \mu_N; \mlim) \int_{\xi'}^{\xi''} P_n(\xi) d\mlim(\xi)
\end{equation}
when $|\xi' - \xi''|$ is small. For $\mlim = \Wig$, a regularisation procedure
suggested in \cite{meband2}  allows to establish (\ref{eq:tomes})
(and even to determine the subleading asymptotic terms) in the cases when (\ref{eq:etcond})
is violated due to divergent contribution to 
\[ \tils(n; \mu_N; \Wig) = \int_{-\infty}^\infty U_n(\xi) d\mu_N(\xi)\]
 coming from the neighbourhood of $\xi = \pm 1$.
It would be interesting to find a way to consider even
shorter scales $\epsilon_N$, for which the limit of 
$\epsilon_N^{-1} R_{\epsilon_N}^{\xi_0} [\mu_N]$  
is distinct from that of $\epsilon_N^{-1} R_{\epsilon_N}^{\xi_0} [\mlim]$. In the random
matrix applications, such a method would allow to study
the local eigenvalue statistics in the bulk of the spectrum via modified moments (in particular,
in problems where alternative methods are not currently available).

\section{Spectral edges of random matrices}\label{s:edges}

\subsection{Wigner matrices}\label{s:edges.Wig}

The application of the moment method to local eigenvalues statistics originates in the work of
Soshnikov \cite{Sosh} on universality for Wigner
matrices. Let us recall the result of \cite{Sosh},
after some preliminaries.

As before, we consider a sequence $(H^{(N)})_{N \geq 1}$ of Wigner matrices, i.e.\ random Hermitian 
matrices such that the diagonal entries
of every $H^{(N)}$ are sampled from a probability
distribution $\mathcal{L}_\text{diag}$ satisfying 
(\ref{eq:diag}), and the off-diagonal entries are 
sampled from a probability distribution $\mathcal{L}_\text{off-diag}$ satisfying  (\ref{eq:offdiag});  the eigenvalues of $H^{(N)}$ are denoted
\[ \xi_1^{(N)} \geq \xi_2^{(N)} \geq \cdots
\geq \xi_N^{(N)}~.\]

Consider the random
point process (i.e.\ a random collection or points, or, equivalently, a random integer-valued
measure)
\begin{equation}\label{eq:edgeresc}
\Lambda^{(N)} = \sum_{j=1}^N \delta\left(
\lambda - N^{1/6} \left[ \xi_j^{(N)} - 2\sqrt{N}\right] \right)
\end{equation}
(the scaling is natural in view of the square-root singularity of $\Wig$ at $1$).

Two special cases, the Gaussian Orthogonal Ensemble (GOE), and the Gaussian Unitary Ensemble
(GUE) [as well as the Gaussian Symplectic Ensemble (GSE, not discussed here)], enjoy an invariance property which allows to apply the method of orthogonal polynomials (see Mehta \cite{Mehta}). The limits of $\Lambda^{(N)}$ for GOE and GUE, called
the Airy$_1$ ($\mathfrak{Ai}_1$) and the Airy$_2$ ($\mathfrak{Ai}_2$) point processes, respectively, were found by Tracy and Widom \cites{TW1,TW2} and Forrester \cite{For1}. The correlation functions, which are (by definition) the densities
\[ \rho_{\beta,k}(\lambda_1, \cdots, \lambda_k)
= \frac{d}{d \, \mes_k} \mathbb{E} \mathfrak{Ai}_\beta^{\otimes k} |_{\lambda_1 < \cdots < \lambda_k}\]
of the off-diagonal parts of the moment measures $\mathbb{E} \mathfrak{Ai}_\beta^{\otimes k}$,
are expressed via determinants involving the Airy function $\Ai$:
\begin{align}
\label{eq:rho2}&\rho_{2,k}(\lambda_1, \cdots, \lambda_k) = 
\det_{k \times k}(A(\lambda_p, \lambda_r))_{p,r=1}^k~,\\
\label{eq:rho1}&\rho_{1,k}(\lambda_1, \cdots, \lambda_k) =
 \sqrt{\det_{2k \times 2k}(A_1(\lambda_p, \lambda_r))_{p,r=1}^k}~,
\end{align}
where
\begin{align*} 
A(\lambda, \lambda') = \int_0^\infty \Ai(\lambda + u) \Ai(\lambda' + u) du~, \,\,
A_1(\lambda, \lambda') = \left(\begin{array}{cc}
A(\lambda, \lambda') &  DA(\lambda, \lambda')\\
 JA(\lambda, \lambda')&A(\lambda, \lambda')
\end{array}
\right)~, \\
DA(\lambda, \lambda') = \frac{\partial}{\partial \lambda'} A(\lambda, \lambda')~, \,\,
JA(\lambda, \lambda') = - \int_{\lambda}^\infty A(\lambda'', \lambda') d\lambda'' - \frac{1}{2} \operatorname{sign}(\lambda-\lambda')~.
\end{align*}

\begin{theorem*}[Soshnikov \cite{Sosh}]
Let $(H^{(N)})_{N\geq1}$ be a sequence of Wigner matrices satisfying the additional assumptions
\begin{align}
\label{eq:sym}&H(u,v) \overset{\text{distr}}{=} -H(u,v)~; &\text{(symmetry)}&\\
\label{eq:subg}&\mathbb{E} |H(u,v)|^{2k} \leq (C k)^k &\text{(subgaussian tails)}&
\end{align} 
on $\mathcal{L}_\text{diag}$ and $\mathcal{L}_\text{off-diag}$. If $\mathcal{L}_\text{off-diag}$
is supported on the real line,
 the point processes $\Lambda^{(N)}$ converge (in the topology of Remark~\ref{rem:top}) to 
 $\mathfrak{Ai}_1$; otherwise, $\Lambda^{(N)} \to \mathfrak{Ai}_2$.
\end{theorem*}

\begin{remark} Lee and Yin \cite{LY} have shown that the 
theorem remains valid if (\ref{eq:sym}) and (\ref{eq:subg}) are replaced with
the assumption
\begin{equation}
\lim_{R \to \infty} R^4 \, \mathbb{P} \left\{ |H(1,2)|  \geq R \right\} = 0~,
\end{equation}
which they have shown to be necessary and sufficient. Their argument makes use of the methods
developed in the works of Erd\H{o}s, Bourgade, Knowles, Schlein, Yau, and Yin on universality in the bulk for Wigner matrices, cf.\  Erd\H{o}s \cite{Erd}.
\end{remark}

\begin{remark}
The work of Soshnikov was followed by numerous other applications of the moment method to
local eigenvalue statistics 
in random matrix theory (see Soshnikov \cite{Sosh2}, P\'ech\'e \cite{Peche}) as well as outside it (see Okounkov \cite{Ok}).
\end{remark}

The strategy of \cite{Sosh} is to compute the asymptotics of moments and to show, using a version of 
Proposition~\ref{p:lapl}, that
the limit of $\Lambda^{(N)}$ exists and does not depend on the distribution of the entries. Thus
the theorem is reduced to its special case 
appertaining to  the Gaussian invariant ensembles.

\subsubsection{An argument based on modified moments}

In \cite{FS1}, modified moments were used to re-prove Soshnikov's theorem quoted above (the method was also
applied to sample covariance matrices, to re-prove the results of Soshnikov \cite{Sosh2} and P\'ech\'e \cite{Peche} on the largest eigenvalues, and to prove a new result on the smallest ones). Let us outline the argument of \cite{FS1} (incorporating modifications from \cite{meband1}), which serves as the basis for the extensions described later in this section.

Let us assume that the diagonal entries $H(u, u)$ are identically zero, and that the off-diagonal entries $H(u, v)$ are randomly chosen signs $\pm 1$. Then Proposition~\ref{prop:cl} identifies 
\begin{equation}\label{eq:eprod}
\mathbb{E} \prod_{r=1}^k \tr P_{n_r}^{(N-1)} \left( \frac{H^{(N)}}{2\sqrt{N-2}} \right)
\end{equation}
as $(N-2)^{-\sum n_r/2}$ times the number of $k$-tuples of closed non-backtracking paths in the complete graph on $N$ vertices, in which
every edge is traversed an even number of times (in total). Such $k$-tuples are divided in topological equivalence classes ($k$-diagrams of Section~\ref{s:edge:constr} below). For $n_r \asymp N^{1/3}$, the contribution of 
every equivalence class can be asymptotically evaluated.

In the regime $n_r \asymp N^{1/3}$,  (\ref{eq:eprod}) captures the asymptotics of the transform (\ref{eq:transf}) of the
moment measures of $\Lambda^{(N)}$. The combinatorial classification yields a convergent series
for this transform. This allows to describe the
vague limit of $\Lambda^{(N)}$. A more general argument making use of an extension of Proposition~\ref{prop:cl} allows to show that the same
limit appears for any sequence of matrices satisfying the assumptions of Soshnikov's theorem (with $ \mathcal{L}_\text{off-diag}$  supported on the real line) in particular, for the GOE for which the answer
is already identified as the Airy$_1$ point process $\mathfrak{Ai}_1$.

\vspace{1mm}\noindent
In the remainder of this section we describe 
(without proofs) two results which may be seen 
as generalisations of \cite{Sosh}.

\subsection{Wigner processes}\label{s:edges.wigproc}

Instead of a single Wigner matrix $H^{(N)}$, let us consider a family $H^{(N)}(\xx) = (H(\xx; u, v))_{1 \leq u \leq v \leq N}$ of Wigner matrices depending on a parameter $\xx \in \XX$; then we are interested in the  eigenvalues 
\[ \xi_1^{(N)}(\xx) \geq \xi_2^{(N)}(\xx)
\geq \cdots \geq \xi_N^{(N)}(\xx)\]  
as a random process on $\XX$.

Let us assume that 
$(\xx \mapsto H(\xx; u, u))_{1 \leq u \leq N}$ are independent copies of a random process $\diag: \XX \to \mathbb{R}$,
\[ \mathbb{E} \, \diag(\xx) = 0~, \quad
\mathbb{E} \, \diag(\xx) ^2  < \infty~, \]
and that $(\xx \mapsto H(\xx; u, v))_{1 \leq u < v \leq N}$ are independent copies of $\offdiag: \XX \to \mathbb{C}$,
\[  \mathbb{E} \, \offdiag(\xx) = 0~, \quad
\mathbb{E} \, |\offdiag(\xx)| ^2  =1~.\]
The process $\offdiag(\xx)$ equips $\XX$ with the $L_2$ metric 
\[ \metr(\xx, \xx') = \sqrt{\frac{1}{2} \mathbb{E} |\offdiag(\xx) - \offdiag(\xx')|^2}~. \]
The local properties of the eigenvalues rescaled about $\xx_0 \in \XX$ depend on the behaviour of $\metr$ near $\xx_0$, which may be captured by the tangent cone  $T_{\xx_0}\XX$ to $\XX$ at $\xx_0$ (the tangent cone to a metric space was introduced by
Gromov \cite{Gromov}*{Section 7}). 

Moment-based methods allow to obtain rigorous results at the spectral edges. Here we focus on the special case 
in which $\XX = \mathbb{R}^d$ and 
\begin{equation}\label{eq:ounorm}
\metr(\xx, \xx')^2 = \|\xx - \xx'\|_p + o(\|\xx\|_p + \|\xx'\|_p) \quad (\xx, \xx' \to 0)
\end{equation}
for some $1 \leq p \leq 2$, which includes, for example, the Ornstein--Uhlenbeck sheet. In this case the tangent cone at the origin is the space $\XX_p^d = (\mathbb{R}^d, \sqrt{\| \cdot\|_p})$.

\begin{theorem}\label{thm:wigproc} Let $d \in \mathbb{N}$ and $1 \leq p \leq 2$. Let $\XX = \mathbb{R}^d$, and suppose the processes $\diag(\xx)$ and $\offdiag(\xx)$ have symmetric distribution  (\ref{eq:sym}) and subgaussian tails  (\ref{eq:subg}) at every point $\xx \in \mathbb{R}^d$, and that the covariance of $\offdiag(\xx)$ has the asymptotics (\ref{eq:ounorm}) near the origin. Then the processes
\[ \Lambda^{(N)}(\xx) = \sum_{j = 1}^N \delta\left(\lambda - N^{1/6}\left[ \xi_j^{(N)}(\xx N^{1/3}) - 2 \sqrt{N}  \right]\right)\]
converge (as $N \to \infty$, in the sense of finite-dimensional distributions) to a limiting process
\[ \AD_\beta[\XX_p^d](\xx) = \sum_{j=1}^\infty \delta(\lambda - \lambda_j(\xx)) \quad (\xx \in \mathbb{R}^d)\]
taking values in sequences $\lambda_1(\xx) \geq \lambda_2(\xx) \geq \cdots$, where $\beta = 1$ if $\offdiag(0)$ is real-valued, and $\beta = 2$ ---  otherwise.
\end{theorem}

The level of generality chosen here is motivated in particular by the following result proved in \cite{mecorner}: the process $\AD_\beta[\XX_1^2]$ also describes the edge scaling limit of corners of time-dependent random matrices (for a discussion of
these, see Borodin \cites{B1,B2}).

\begin{proposition} For $\beta \in \{1, 2\}$, the  process $\AD_\beta[\XX_p^d]$ boasts the  properties:
\begin{enumerate}
\item
There exists a modification of $\AD_\beta[\XX_p^d](\xx)$ in which every $\lambda_j(\xx)$ is a continuous function of $\xx \in \mathbb{R}^d$.
\item 
At a fixed $\xx \in \mathbb{R}^d$,  $\AD_\beta[\XX_p^d](\xx)$ is equal in distribution to the Airy$_\beta$ point process $\mathfrak{Ai}_\beta$.
\item The distribution of $\AD_\beta[\XX_p^d](\xx)$
at a $k$-tuple of points $(\xx_q)_{q=1}^k$ in $\mathbb{R}^d$ depends only on $\beta$ and on the distances $\|\xx_q - \xx_r\|_p$ ($1 \leq q < r \leq k$).
\end{enumerate}
\end{proposition}

The last item implies that the distribution of the restriction of $\AD_\beta[\XX_p^d]$ to a geodesic in $\ell_d^p$ does not depend on the choice of geodesic
(and neither on $p$ and $d$), and thus coincides in distribution with $\AD_\beta[\XX^1]$. 

We note that $\AD_2[\XX^1]$ admits a concise determinantal description. Indeed, the $\beta=2$ Dyson Brownian motion satisfies the assumptions of 
Theorem~\ref{thm:wigproc}; thus its edge
scaling limit is described by the process $\AD_2[\XX^1]$. On the other hand, Mac\^edo \cite{Mac}
and Forrester, Nagao, and Honner \cite{FNH} (see further Forrester \cite{Forbook}*{7.1.5}) found this limit directly. This process,
 the moment measures of which are given by determinants, appeared again  in the work of Pr\"ahofer and Spohn \cite{PrSp} on models of random growth (see the lecture notes of Johansson \cite{Johdet} for further limit theorems in which it appears);  Corwin and Hammond \cite{CorHam} studied its properties, and  coined the  term `Airy line ensemble'. Thus the distribution of the restriction of $\AD_2[\XX_p^d]$ to any geodesic in $\ell_p^d$ is given by the Airy line ensemble.

\subsubsection{Construction of the processes $\AD_\beta$}\label{s:edge:constr}

With the exception of the case $\beta = 2$, $d=1$,
the process $\mathfrak{P} = \AD_\beta[\XX_p^d]$ does not seem to be described by  determinantal formul{\ae}. The construction presented here is motivated by the combinatorial arguments of Soshnikov \cite{Sosh} and further by the work of Okounkov \cite{Ok}, and makes use of the
results of \cite{FS1}.

Let  $\mathfrak{P}(\xx)$ ($\xx \in \mathbb{R}^d$) 
be a random process which takes values in point configurations on the line (i.e.\ locally
finite sums of $\delta$-functions). That is, for every $\xx \in \mathbb{R}^d$ the random variable $\mathfrak{P}(\xx)$ is a point process on $\mathbb{R}$. Denote
\begin{equation} \widetilde{\rho}_{\mathfrak{P}, k}
(\xx_1, \cdots, \xx_k) = \mathbb{E} \prod_{r=1}^k \mathfrak{P}(\xx_r) \quad
(\xx_1, \cdots, \xx_k \in \mathbb{R}^d)
\end{equation}
be the moment measures of $\mathfrak{P}$, and
\[
\widetilde{R}_{\mathfrak{P}, k}(\xx_1, \cdots, \xx_k; \alpha_1, \cdots, \alpha_k)= \int \prod_{r=1}^k \frac{\sin \alpha_r \sqrt{-\lambda_r}}{ \sqrt{-\lambda_r}} d\widetilde{\rho}_{\mathfrak{P},k}(\xx_1, \cdots, \xx_k; \lambda_1, \cdots, \lambda_k)
\]
(in our case, the  divergent integral can be regularised, cf.\ \cite{meband1}*{Section~6}; the transform appears from the asymptotics of orthogonal polynomials, cf.\ (\ref{eq:transf})). Then, let 
\begin{equation}
\widetilde{R}_{\mathfrak{P},k}^\# (\bar{\xx}; \alpha)
= \sum_{I \subset \{1, \cdots, k\}}
\widetilde{R}_{\mathfrak{P},\#I}(\bar{\xx}|_I, \alpha|_I) 
\widetilde{R}_{\mathfrak{P},k-\#I}(\bar{\xx}|_{I^c},\alpha|_{I^c})~,
\end{equation}
where $\bar{\xx}=(\xx_1, \cdots, \xx_k) \in (\mathbb{R}^d)^k$, and
$\bar{\xx}|_I = (\xx_r)_{r\in I}$. (The sum over partitions has to do with the contribution of the two spectral edges to the asymptotics.)
We define $\AD_\beta$ via a formula (\ref{eq:adtransf})  for $\widetilde{R}^\#_{\AD_\beta[\XX_p^d],k}$ (uniqueness follows
from the considerations of  \cite{meband1}*{Section~6}). 

\begin{figure}[ht]
\vspace{2.2cm}
\setlength{\unitlength}{1cm}
\begin{pspicture}(0,0)
\psline(.2,1)(2, 1)
\psarcn(2.8,1){.8}{180}{195}
\psset{arrows=-}
\psline(2,.8)(2.1,.92)
\psarcn(2.8,1){.7}{185}{195}
\psline(2,.92)(2.1,.8)
\psset{arrows=->}
\psline(2,.9)(.2,.9) %
\psline(4.2,1)(6,1)
\psarcn(6.8,1){.8}{180}{185}
\psline(6,.9)(5,.9)
\pscurve(5,.9)(5.5,.7)(6,.5)(6.4,.5)
\psarc(6.8,1){.65}{230}{218}
\pscurve(6.27,.6)(5.9, .6)(5,.8)(4.8,.9)(4.2,.9) %
\psline(8.2,1)(9, 1)
\psline(9, 1)(9.2, 1.7)
\psarcn(9.3, 2.05){.35}{255}{270}
\psarcn(9.3, 2.05){.45}{265}{270}
\psline(9.3, 1.7)(9.1, 1)
\psline(9.1, 1)(10, 1)
\psarcn(10.8,1){.8}{180}{195}
\psset{arrows=-}
\psline(10,.8)(10.1,.92)
\psarcn(10.8,1){.7}{185}{195}
\psline(10,.92)(10.1,.8)
\psset{arrows=->}
\psline(10,.9)(8.2,.9)
\end{pspicture}
\caption{Three $1$-diagrams. The diagram on the left corresponds to the projective plane with $s = 1$ (left). The two diagrams
in the centre and on the right correspond to surfaces with $s=2$; the one is the centre is the torus.}\label{fig:diag1}
\vspace{2mm}
\end{figure}
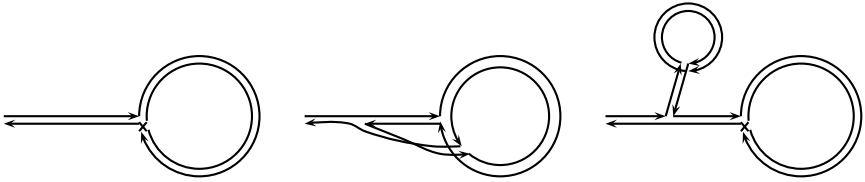

\vspace{2mm}\noindent
Let us consider the collection of $k$-tuples of 
non-backtracking walks for which every edge of
the spanned graph $\mathfrak{g}$ (cf.\ Section~\ref{sub:Wigner}) is traversed 
exactly twice, and every vertex in $\mathfrak{g}$
has degree at most three. Such $k$-tuples can
be divided into topological equivalence classes ($k$-diagrams). For example, for
$k=1$, the paths $1 \, 2 \, 3 \, 4 \,\- 5 \, 6 \, 7 \, 8 \, 4 \, 5 \, 6 \, 7 \, 8 \, 4 \, 3 \, 2 \, 1$ on Figure~\ref{fig:forest} (right) and the path
$8 \, 7 \, 2 \, 3 \, 9 \, 2 \, 3 \, 9 \, 2 \, 7 \, 8$
belong to the same equivalence class, schematically
depicted on Figure~\ref{fig:diag1} (left). The formal 
definition is given in \cites{FS1,meband1}.

\begin{figure}[ht]
\vspace{2.5cm}
\setlength{\unitlength}{1cm}
\begin{pspicture}(-1,0)
\psline(.2,1)(1, 1)
\psarcn(1.6,1){.6}{180}{187}
\psline(1,.9)(.2,.9)
\psline(3,1)(2.05, 1)
\psarc(1.6,1){.5}{0}{352}
\psline(2.05, .9)(3,.9)
\psline(4,0)(4, 1)
\psarcn(4,1.5){.5}{270}{278}
\psset{arrows=-}
\psline(4.1,1)(4,1.1)
\psarcn(4,1.5){.4}{268}{282}
\psline(4.1,1.1)(4,1)
\psset{arrows=->}
\psline(4.1,1)(4.1,0)
\psline(4.7,2)(4.7, 1)
\psarc(4.7,.5){.5}{90}{82}
\psset{arrows=-}
\psline(4.8,1)(4.7,.9)
\psarc(4.7,.5){.4}{90}{82}
\psline(4.8,.9)(4.7,1)
\psset{arrows=->}
\psline(4.8,1)(4.8,1.5)
\psline(4.8,1.5)(5.8, 1.5)
\psarc(6.4,1.5){.6}{180}{173}
\psset{arrows=-}
\psline(5.8, 1.6)(5.9, 1.5)
\psarc(6.4,1.5){.5}{180}{170}
\psline(5.9,1.6)(4.8, 1.6)
\psset{arrows=->}
\psline(4.8,1.6)(4.8,2)
\psline(7.7,1)(8.5, 1)
\rput(8, 1.2){\tiny I}
\psline(8.5, 1)(8.7, 1.7)
\rput(8.4, 1.35){\tiny II}
\psarcn(8.8, 2.05){.35}{255}{270}
\psarcn(8.8, 2.05){.45}{265}{270}
\rput(8.2, 2.35){\tiny III}
\psline(8.8, 1.7)(8.6, 1)
\psline(8.6, 1)(9.5, 1)
\rput(9.1, 1.2){\tiny IV}
\psarcn(10.1,1){.6}{180}{190}
\rput(10.1, 1.8){\tiny V}
\rput(10.1, 0.7){\tiny VI}
\psline(9.5,.9)(7.7,.9)
\psline(11.1, 2)(10.45, 1.35)
\psarc(10.1,1){.5}{45}{37}
\psline(10.5, 1.25)(11.17, 1.93)
\rput(10.65, 1.85){\tiny VII}
\end{pspicture}
\caption{Three $2$-diagrams: $s = 2$ (left),
    $s = 3$ (centre, right). The leftmost diagram, corresponding to a 
sphere glued from two disks, is often responsible for fluctuations of linear eigenvalue statistics on global and mesoscopic scales.}\label{fig:diag2}
\end{figure}
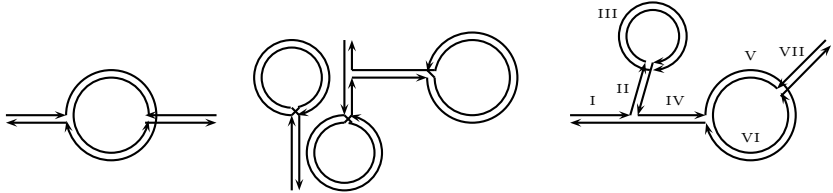

Different $k$-diagrams correspond to homo\-to\-pi\-cal\-ly
distinct ways to glue $k$ disks with a marked point
on the boundary. The result of such a gluing is
a two-dimensio\-nal manifold. Thus to every $k$-diagram one
can associate a number $s$, which is related to the
Euler characteristic $\chi$ of the manifold by the formula $s = 2 k - \chi$; for $k=1$ the number $s$ is  the non-oriented genus. 
 The multi-graph associated to a diagram with a certain value of $s$ (see  
Figures~\ref{fig:diag1} and \ref{fig:diag2}) has $2s$ vertices and $3s-k$
edges. The number $D_k(s)$ of $k$-diagrams with a given value of $s$ satisfies the estimates (\cite{FS1}*{Proposition~II.3.3})
\begin{equation} \frac{(s/C)^{s+k-1}}{(k-1)!} \leq  D_k(s) \leq \frac{(Cs)^{s+k-1}}{(k-1)!}~;\end{equation}
the upper bound guarantees that the series (\ref{eq:adtransf}) which we derive below converges.

Next, we associate to a $k$-diagram $\mathcal{D}$ and to $\alpha \in (0, \infty)^k$ a 
$(3s-2k)$-dimensional polytope $\Delta_{\mathcal{D}}(\alpha)$ in $\mathbb{R}^{3s-k}$,
as follows. The variables $w(e)$ are labeled by the edges $e$ of $\mathcal{D}$; the polytope is defined by the inequalities 
\[\left\{\begin{aligned}
&w(e) \geq 0 &(e \in \mathrm{Edges}(\mathcal{D}))\\
&\sum_{e} c_r(e) w(e) = \alpha_r &(1 \leq r \leq k)
\end{aligned}\right.\]
where $c_r(e) \in \{0,1,2\}$ is the number of times the edge $e$ is traversed by the $p$-th path in the diagram. For example, the polytope associated with
the rightmost $2$-diagram of Figure~\ref{fig:diag2} is
given by
\[ \left\{ \begin{aligned}
&w(\text{I}), w(\text{II}), \cdots, w(\text{VII}) \geq 0\\
&2w(\text{I}) + 2w(\text{II}) + 2w(\text{III}) + 2w(\text{IV}) + w(\text{V}) + w(\text{VI}) = \alpha_1 \\
& w(\text{V}) + w(\text{VI}) + 2 w(\text{VII}) = \alpha_2~.
\end{aligned} \right. \]

Let $\mathfrak{D}_1(k)$ be the collection of all $k$-diagrams, and let  $\mathfrak{D}_2(k) \subset  \mathfrak{D}_1(k)$ be the sub-collection of diagrams in which every edge is traversed once in one direction and once in another one (such as Figure~\ref{fig:diag1}, centre, and Figure~\ref{fig:diag2}, left; these diagrams correspond to gluings preserving orientation). Now we can finally write the series for $\widetilde{R}^\#$:

\begin{multline}\label{eq:adtransf}
\widetilde{R}_{\AD_\beta,k}^\# (\bar{\xx}; \alpha) = \\ \sum_{\mathcal{D} \in \mathfrak{D}_\beta} \int_{\Delta_{\mathcal{D}}(\alpha)} \, \exp\left\{- \sum_{e \in \mathrm{Edges}(\mathcal{D})}  \| \xx_{r_+(e)} - \xx_{r_-(e)}\|_p w(e) \right\} \, d\,\mes_{3s-2k}(w)~,
\end{multline}
where  $k \geq r_+(e) \geq r_-(e) \geq 1$ are the indices of the two paths traversing $e$ in $\mathcal{D}$. 

For example, when $\XX$ is a singleton, all the terms in the exponent vanish, and (\ref{eq:adtransf}) yields an expression for the Airy point process in terms of volumes of the polytopes $\Delta_{\mathcal{D}}(\alpha)$,
which may be compared to the one given by Okounkov \cite{Ok}*{\S2.5.4}.

\subsection{Band matrices}\label{s:edges.band}

In this section, we discuss an extension of Soshnikov's theorem to a class of matrices of the form
considered in Section~\ref{sub:Wigner}. First, we recall a conjecture, based on the Thouless 
criterion \cite{Th}. Then we discuss a particular case, the spectral
edges of random band matrices, in which the conjecture can be  proved. Finally, we comment on
mesoscopic scales.

\subsubsection{Thouless criterion}\label{s:edges.Th} 

 The Thouless criterion \cite{Th}, originally introduced in the context of Anderson localisation, can be applied to predict the behaviour of local eigenvalue statistics; cf.\ Fyodorov
and Mirlin \cites{FM,FM2}. Consider a sequence of matrices $H^{(N)}$ associated with a sequence of 
graphs $G_N = (V_N, E_N)$ as in Section~\ref{sub:Wigner}. Then the measures
\[ \mu_N = \frac{1}{\# V_N} \sum_{j=1}^N \delta\left(\xi - \frac{\xi_j^{(N)}}{2\sqrt{2W_N}}\right) \] 
converge to the semi-circle measure $\Wig$.

Let $\xi_0 \in \mathbb{R}$, and let $\eta_N > 0$ be chosen so that 
the sequence of (random) measures 
\[ \# V_N R_{\eta_N}^{\xi_0}[\mathbb{E} \mu_N] 
= \mathbb{E} \sum_{j=1}^{\# V_N} \delta\left( \lambda - \frac{1}{\eta_N} \left[  \frac{\xi_j^{(N)}}{2 \sqrt{2 W_N-1} } - \xi_0 \right]\right) \] 
will have a non-trivial vague
limit (cf.\ Section~\ref{s:resc}). Thus chosen, $\eta_N$ measures the mean spacing between 
eigenvalues, whereas
\[ (\mathrm{spacing/DOS})\,(\xi_0) = \eta_N^2 \#V_N \]
measures the mean spacing in units of the density of states. Let us compare the inverse of this quantity with
the mixing time $T^{\mathrm{mix}}$ of the random walk on $G_N$.\footnote{equivalently, the ratio of the mixing time and the density of states, which is interpreted as the energy-dependent 
mixing time, is compared to the usual inverse eigenvalue spacing $\eta_N$.}
In many cases the following seems to be correct: the eigenvalue statistics of $H^{(N)}$ near $\xi_0$
are described by random matrix theory if and only if
\begin{equation}\label{eq:thouless}
T^\text{mix}(G_N) \ll \frac{1}{(\mathrm{spacing/DOS})(1)}~.
\end{equation}
This interpretation of the Thouless criterion is based
on the assumption that the semi-classical approximation is valid up to the scales governing
the local eigenvalue statistics; we refer to the reviews
of Spencer \cites{Sp-banded,Sp-And} for a discussion
of various aspects of Thouless scaling and its mathematical justification, and to the work of 
Spencer and Wang \cite{Wang} for some rigorous results. Here we focus our attention
on the particular case of

\subsubsection{Random band matrices}
Denote
$\|u - v\|_N = \min_{\ell \in \mathbb{Z}} |u - v - \ell N|$.
A (one-dimensional) random band matrix of bandwidth $W$ is for us a random Hermitian $N \times N$ matrix $H^{(N)} = (H(u, v))_{1 \leq u,v \leq N}$
such that
\begin{equation} \label{eq:band}\begin{cases}
H(u, v) = 0~, &\|u-v\|_N > W~,\\
H(u, v) \sim \mathcal{L}_\text{off-diag}~,
&1 \leq \|u-v\|_N \leq W ~, \\
H(u, u) \sim \mathcal{L}_\text{diag}~,
\end{cases}\end{equation}
where $\mathcal{L}_\text{diag}$ and $\mathcal{L}_\text{off-diag}$ satisfy the normalisation conditions  (\ref{eq:diag}) and (\ref{eq:offdiag}), respectively. In the setting
of Section~\ref{sub:Wigner}, it corresponds to the
graph $G_N = (V_N, E_N)$,
\begin{equation}\label{eq:defgr} V_N = \{1, \cdots, N\}~, \quad
(u, v) \in E_N \iff 1 \leq \|u-v\|_N \leq W_N~,
\end{equation}
More general band matrices are discussed, for example, in 
\cites{KMP,Sp-banded,EK2}.

For $-1 < \xi_0 < 1$ (the bulk of the spectrum), 
\[(\mathrm{spacing/DOS})\,(\xi_0) \asymp \frac{1}{N}~,
\quad T^{\mathrm{mix}} \asymp \frac{N^2}{W^2}~, \]
therefore the criterion (\ref{eq:thouless}) suggests
the following: the eigenvalue statistics of $H^{(N)}$ near $\xi_0$ are described by 
random matrix theory if and only if $W \gg \sqrt{N}$. This prediction is supported
by the detailed super-symmetric analysis performed by Fyodorov and Mirlin \cites{FM,FM2}.
Mathematical justification remains a major challenge, cf.\ Spencer \cites{Sp-banded,Sp-SUSY}
and references therein.

\subsubsection{Spectral edges}

The (modified) moment method allows to confirm  the criterion (\ref{eq:thouless})
at the spectral edges of random band matrices.  

\begin{theorem}[cf. \cite{meband1}*{Theorem~1.1}]\label{thm:band} Let $(H^{(N)})_{N \geq1}$ be a sequence of random band matrices satisfying the unimodality assumptions (\ref{eq:unimod}). If the 
bandwidth $W_N$ of $H^{(N)}$ satisfies
\begin{equation}\label{eq:cond56}
 \lim_{N \to \infty} \frac{W_N}{N^{5/6}} = \infty~, 
\end{equation}
then
\[ \sum_{j=1}^N \delta \left( \lambda - \frac{N^{2/3}}{\sqrt{2W_N}} \left[ \xi_j^{(N)} - 2 \sqrt{2W_N} \right] \right) \to \mathfrak{Ai}_\beta~,\]
where $\beta = 1$ if $\supp \mathcal{L}_{\text{off-diag}} \subset \mathbb{R}$, and $\beta = 2$ otherwise.
\end{theorem}
The threshold $N^{5/6}$ in (\ref{eq:cond56}) is sharp, see \cite{meband1}*{Theorem~1.2}. The same
\cite{meband1}*{Theorem~1.2} implies that
\[ \eta_N \asymp \min(W_N^{2/5}N^{-1}, N^{-2/3})~, \]
therefore
\[ (\mathrm{spacing/DOS})(1) \asymp  \min(W_N^{4/5}N^{-1}, N^{-1/3})~,\] 
and (\ref{eq:cond56}) is consistent with (\ref{eq:thouless}).

The unimodality
conditions (\ref{eq:unimod}) simplify the analysis 
(cf.\ Proposition~\ref{prop:cl}); we expect that they can be relaxed using the methods of \cite{FS1}*{Part~III} and \cite{EK2}.

\subsubsection{Mesoscopic scales}\label{sub:EK} 
On mesoscopic scales $1 \gg \epsilon_N \gg 1/\#V_N$, the following counterpart of the Thouless criterion
goes back to the (physical) work of Altshuler and Shklovskii \cite{AS}. Let $\eta_N$ be such that  the sequence
$(\epsilon_N^{-1} R_{\eta_N}[\mathbb{E} \mu_N])_N$ has a non-trivial vague limit. If
\begin{equation}\label{eq:meso}
\epsilon_N \eta_N^{-2} \gg T^\text{mix}~,
\end{equation}
 the fluctuations of linear eigenvalue statistics should be described by 
a log-correlated Gaussian field, whereas when (\ref{eq:meso}) is violated, one expects a more
regular field depending on the geometry of the underlying lattice. We refer to the works of Fyodorov, 
Le Doussal, and Rosso \cite{FLR} and of Fyodorov and Keating \cite{FK} for a discussion
of the significance of log-correlated fields within
and outside random matrix theory, and to the work
of Fyodorov, Khoruzhenko, and Simm \cite{FKS} for results 
pertaining to the Gaussian Unitary Ensemble.

Erd\H{o}s and Knowles proved a series of  results pertaining to mesoscopic statistics for a wide
class of $d$-dimensional band matrices. In the works \cites{EK1,EK2}, they developed a 
moment-based approach which allowed them to control the quantum dynamics associated for time scales 
$t \leq W_N^{d/3 - \delta}$.  In \cites{EK3,EK4}, they gave mathematical justification to the  criterion 
(\ref{eq:meso}) in the range $\epsilon_N \geq W_N^{- d/3 + \delta}$. It would be interesting to extend the
results of \cites{EK1,EK2} and \cites{EK3,EK4} to the full mesoscopic range.

\section{Some further questions}

\noindent{\bf Other limiting measures} The spectral
measures in this article converge to the semicircle distribution $\mlim = \Wig$. The modified moment method 
described here has been also applied to the Kesten--McKay measure
(the orthogonality measure for $P_n^{(\kappa)}$), the Godsil--Mohar measure
(its bipartite analogue), and the Marchenko--Pastur
measure (the infinite connectivity limit of the Godsil--Mohar measure); see e.g.\ \cites{meop,mesparse}. It would be interesting to
adapt the method to situations in which the recurrent relation has less explicit form.

\vspace{2mm}\noindent{\bf $\beta$-ensembles} The (convincing, although so far unrigorous) ghost and shadows formalism introduced by Edelman \cite{Edelman}
strong\-ly suggests that the construction (\ref{eq:adtransf}) should have an extension to
general $\beta>0$. See Forrester \cite{Forbeta} and
\cite{Forbook} for background on $\beta$-ensembles,
and Borodin and Gorin \cite{BorGor} for a recent result
pertaining to the spectral statistics of submatrices of $\beta$-Jacobi random matrices.

\vspace{2mm}\noindent 
{\bf Time-dependent invariant ensembles} It seems plausible that, for general (non-Gaussian) invariant ensembles undergoing  Dyson-type evolution, the spectral statistics near a soft edge should be described by the processes $\AD_\beta$ of Section~\ref{s:edge:constr}. Currently, there seem to be no proved results of this form (even for the case $\beta=2$ in which determinantal formal{\ae} for finite
matrix size are given by the Eynard--Mehta theorem \cite{Mehta}*{Chapter 23}).

\vspace{2mm}\noindent{\bf Beyond random matrices} Motivated by
the proof of the Baik--Deift--Johans\-son conjecture given by Okounkov \cite{Ok}, one may look
for the appearance of (\ref{eq:adtransf}) outside random matrix theory, particularly,
in the context of random growth models,  for a discussion of the subtle connection between which and random matrix theory we refer to the lecture notes of Ferrari \cite{Fer}.

\vspace{2mm}\noindent{\bf Bulk of the spectrum} 
We are not aware of any derivation (rigorous or not) of the local eigenvalue statistics
in the bulk of the spectrum using any version of the moment method. 
Even for the test case of the Gaussian Unitary Ensemble (tractable by other means), perturbative methods such as Chebyshev expansions have not been of use  beyond the scales $\epsilon_N \gg N^{-1 + \delta}$.
For random band matrices the expansion 
(\ref{eq:L2exp}) has been only regularised for $\epsilon_N \gg W^{-1 + \delta}$ (see \cite{meband2}). 

\vspace{-1mm}
\paragraph{Acknowledgment} 
It is a great pleasure to thank Leonid Pastur, who introduced me to the theory of random matrices and
encouraged me to work on local eigenvalue statistics,
and  Tom Spencer, who spent an immeasurable amount of time to share with me his knowledge, intuition, and taste in various parts of mathematical physics.

 A significant part of the results stated here are based on the joint work \cites{FS1,FS2} with Ohad Feldheim. Yan Fyodorov explained me the relation between (\ref{eq:cond56}) and Thouless scaling.
The results in Section~\ref{s:edges.wigproc} are motivated by the talks given by Alexei Borodin
at the IAS, and by the subsequent discussions with him and with Vadim Gorin. Yan Fyodorov, Vadim Gorin, Antti Knowles, Alon Nishry, Mira Shamis, Misha Sodin, and Ofer Zeitouni kindly commented on a preliminary version of this text. I thank them  very much.

\end{document}